\documentclass[bsl]{asl} 
\usepackage{url}

\title{On the Ja\'{s}kowski models for intuitionistic propositional logic}
\author{R.D. Arthan}
\revauthor{Arthan, R.D.}
\address{Queen Mary University of London\\
Mile End Road\\
London\\
E1 4NS, UK.}

\email{r.arthan@qmul.ac.uk}

\usepackage{amssymb}

\usepackage{prooftree}

\newtheorem{Theorem}{Theorem}
\newtheorem{Lemma}[Theorem]{Lemma}
\newtheorem{Corollary}[Theorem]{Corollary}
\def\Proof{\begin{proof}}
\def\Done{\end{proof}}

\def\Lab#1{{\em(#1)}}

\def\Func#1{{\mathsf{#1}}}

\def\Bar{\mathrel{|}}

\def\Var{\Func{Var}}

\def\Reg#1#2{\begin{array}{|c|}\hline #1 \\ \hline #2 \\\hline \end{array}}

\def\Hide#1{\relax}

\def\AA{\mathbb{A}}
\def\BB{\mathbb{B}}

\def\cL{\mathcal{L}}

\def\cV{\mathcal{V}}

\def\VH{\mathbf{H}}

\def\VJ{\mathbf{J}}

\def\frJ{\mathfrak{J}}

\def\Imp{\Rightarrow}
\def\Pmi{\Leftarrow}
\def\Iff{\Leftrightarrow}
\def\And{\land}
\def\Or{\lor}
\def\Not{\lnot}

\def\SDef{\mathrel{{:}{\equiv}}}
\def\Pr{\vdash}

\def\Diff{\setminus}

\def\Res{\rightarrow}
\def\Meet{\sqcap}
\def\Join{\sqcup}
\def\aF{\Func{f}}
\def\aT{\Func{t}}

\def\IPL{\mathrm{\mathbf{IPL}}}
\def\IPLpr{\IPL\Pr}
\def\Var{\cV}
\def\lF{\bot}
\def\lT{\top}
\def\Star{{*}}
\def\Jas{\VJ}

\def\Eqp{\mathrel{{\dashv}{\vdash}}}

\def\Jaskowski{Ja\'{s}kowski}

\begin{document}

\begin{abstract}

In the 1930s, Stanislaw {\Jaskowski} discovered an
interesting sequence $\frJ_0, \frJ_1, \ldots$ of what he called ``matrices'' and
that today we would think of as finite Heyting Algebras.  He gave a very
brief sketch of a proof that if a propositional formula holds in every $\frJ_i$
then it is provable in intuitionistic propositional logic ($\IPL$).  The sketch just
describes a certain normal form for propositional formulas and gives a very
terse outline of an inductive proof that an unprovable formula in
the normal form can be refuted in one of the $\frJ_k$.  Unfortunately, it is
far from clear how to recover a complete proof from this sketch.

In the early 1950s, Gene F. Rose published a detailed proof of
{\Jaskowski}'s result, still using the notion of matrix rather than Heyting
algebra, based on a normal form that is more restrictive than the one that
{\Jaskowski} proposed.  However, Rose's paper refers to his thesis
for additional details, particularly concerning the normal form.

This note gives a proof of \Jaskowski's result using modern terminology and a
normal form more like \Jaskowski's.  We also prove a semantic property of the
normal form enabling us to give a novel proof of completeness of $\IPL$ for the
Heyting algebra semantics. We outline a decision procedure for $\IPL$ based on
the proof of {\Jaskowski}'s result and illustrate it in action on some simple
examples.

\end{abstract}

\maketitle

Let $\VH = (H, \aF, \aT, \Meet, \Join, \Res)$ be a Heyting algebra. We will
define a new Heyting algebra $\Gamma(\VH)$ by adding a co-atom, i.e., a new
element $*$ such that $x < * < \aT$ for $x \in H \Diff \{\aT\}$.
$\Gamma(\VH)$ will extend $\VH$ as a $(\aF, \aT, \Meet, \Res)$-algebra and the
join in $\Gamma(\VH)$ will agree with the join in $\VH$ wherever possible.
Thus, we choose some object $\Star = \Star_H$ that is not an element of $H$ and
let $\Gamma(\VH) = (H \cup \{\Star\}, \aF, \aT, \Meet, \Join, \Res)$, where the
operations $\Meet$, $\Join$ and $\Res$ are derived from those of $\VH$ as shown
in the operation tables below, in which $x$ and $y$ range over $H \Diff
\{\aT\}$ and where $\alpha : H \to (H \Diff \{\aT\}) \cup \{\Star\}$ satisfies
$\alpha(x) = x$ for $x \neq \aT$ and $\alpha(\aT) = *$.
%\[
%\alpha(x) =
%   \left\{
%      \begin{array}{ll}
%         x & \mbox{if $x \neq \aT$} \\
%         \Star_H & \mbox{if $x = \aT$.}
%      \end{array}
%   \right.
%\]

\[
\begin{array}{ccc}
\begin{array}{c|ccc}
\Meet   & y         & \Star & \aT   \\ \hline
  x     & x \Meet y &   x   & x     \\
\Star   & y         & \Star & \Star \\
 \aT    & y         & \Star & \aT
\end{array}
&
\begin{array}{c|ccc}
\Join   & y         & \Star & \aT   \\ \hline
  x     & \alpha(x \Join y) & \Star & \aT   \\
\Star   & \Star     & \Star & \aT   \\
 \aT    & \aT       & \aT   & \aT
\end{array}
&
\begin{array}{c|ccc}
\Res    & y         & \Star & \aT   \\ \hline
  x     & x \Res y  & \aT   & \aT   \\
\Star   & y         & \aT   & \aT   \\
 \aT    & y         & \Star & \aT
\end{array}
\end{array}
\]
\noindent
%$\Gamma(\VH)$ is the unique extension of $\VH$ to $H \cup \{\Star\}$ as a
%$(\aF, \aT, \Meet, \Res)$-algebra such that {\em(i)} $*$ is a co-atom (i.e., $x
%< \Star < \aT$ for $x \in H \Diff \{\aT\}$) and {\em(ii)} the inclusion of $H$
%in $H \cup \{\Star\}$ preserves joins except where the requirement that $*$ be
%a co-atom forces $p \Join q$ to be $*$ rather than $\aT$.

Let $\BB$ be the two-element Heyting algebra and, as usual, let us write
$\VH^i$ for the $i$-fold power of a Heyting algebra $\VH$.
Then define a sequence $\Jas_0, \Jas_1, \ldots$ of finite Heyting algebras as
follows:
\begin{align*}
\Jas_0 &= \BB \\
\Jas_{k+1} &= \Gamma(\Jas_k^{k+1})
\end{align*}

We take the language $\cL$  of intuitionistic propositional logic, $\IPL$,
to be constructed from a set $\Var = \{P_1, P_2, \ldots\}$ of variables, the
constants $\lF$, $\lT$, and the binary connectives $\And$, $\Or$ and $\Imp$.
We do not take negation as primitive: $\Not A$ is an abbreviation for $A\Imp
\lF$.  The metavariables $A, B, \ldots, M$ (possibly with subscripts) range
over formulas. $E$ and $F$ are reserved for formulas that are either variables
or $\lF$.  $P, Q, \ldots, Z$ range over variables.  We assume known one of the
many ways of defining the logic of $\IPL$ and write $\IPLpr A$, if $A$ is
provable in $\IPL$.  $\IPL$ has an algebraic semantics in which, given a
Heyting algebra $\VH$ and an interpretation $I : \Var \to H$, we extend $I$ to
a mapping $v_I : \cL \to H$ by interpreting $\lF$, $\lT$, $\And$, $\Or$ and
$\Imp$ as $\aF$, $\aT$, $\Meet$, $\Join$ and $\Res$ respectively.  As usual we
write $I \models A$ if $v_I(A) = \aT$, $\VH \models A$ if $I \models A$ for
every interpretation $I : \Var \to H$ and $\models A$ if $\VH \models A$ for
every Heyting algebra $\VH$.  We assume known the fact that $\IPL$ is sound
with respect to this semantics in the sense that, if $\IPL \Pr A$, then
$\models A$.  The converse statement, i.e., the completeness of $\IPL$ with
respect to the semantics is well-known, but we do not use it: in fact we will
give an alternative to the usual proofs.

We write $A \Iff B$ for $(A \Imp B) \And (B \Imp A)$
and $A[B/X]$ for the result of substituting $B$ for each occurrence
of $X$ in $A$. We have the following substitution lemma:

\begin{Lemma}[substitution]\label{lma:subst}
For any formulas $A$, $B$ and $C$ and any variable $X$ we have:

\begin{tabular}{rl}
\Lab{i} &
if $\IPLpr C$, then $\IPLpr C[A/X]$; \\
\Lab{ii} &
if $\IPLpr A \Iff B$, then $\IPLpr C[A/X] \Iff C[B/X]$; \\
\end{tabular}
\end{Lemma}
\Proof {\Lab{i}} is proved by induction on a proof of $C$.
{\Lab{ii}} is proved by induction on the structure of $C$.
\Done

We say a formula $A$ is {\em reduced} if $\lT$ does not appear in $A$
as the operand of any connective and $\lF$ does not appear in $A$
as the operand of any connective other than as the right-hand operand
of $\Imp$.
Thus the only reduced formula containing $\lT$ is $\lT$ itself,
while $\lF$ is only uwed in a reduced formula to form negations.

\begin{Lemma}\label{lma:red}
Any formula is equivalent to a reduced formula.
\end{Lemma}
\Proof This follows by repeated
use of the substitution lemma and the provable equivalences
$\lT \And A \Iff A$, $\lF \And A \Iff \lF$ etc.
\Done

We define a formula to be {\em basic} if it is reduced and is either a
variable or has one of the forms $P \Imp A$ or $A \Imp P$ where $P$ is a
variable and $A$ contains at most one connective.
Thus a basic formula has one of the following forms\footnote{
We elide brackets using the rules that $\Imp$ is right associative
and that the connectives are listed in increasing order of precedence  
as $\Iff$, $\Imp$, $\Or$, $\And$, $\Not$.
}.
\[
\begin{array}{cccccc}
 P & P \Imp Q &P \Imp Q \And R & P \Imp Q \Or R &
   P \Imp Q \Imp R & P \Imp \Not Q  \\
  & \Not P & P \And Q \Imp R & P \Or Q \Imp R &
   (P \Imp Q) \Imp R & \Not P \Imp Q
\end{array}
\]
Note that if $A$ is basic formula of a form other than $P$,
$(P \Imp Q) \Imp R$ or $\Not P \Imp Q$, then $V_I(A) = \aT$ in any Heyting
algebra under the interpretation $I$ that maps every variable to $\aF$.
Our convention for the metavariables  $E$ and $F$ allows us to write, for example, $(P \Imp E) \Imp R$ as a metanotation for the forms $(P \Imp Q) \Imp R$ and $\Not P \Imp R$.

We say a formula is a {\em basic context} if it is reduced and is a conjunction
of one or more pairwise distinct basic formulas.  We say a formula is {\em
regular} if it is an implication $K \Imp F$ where $K$ is a basic context (and
following our convention $F$ is a variable or $\lF$).

We say $A$ and $B$ are {\em equiprovable} and write $A \Eqp B$ if $\IPLpr A$ iff $\IPLpr B$.

\begin{Lemma}\label{lma:reg}
Every formula $A$ is equiprovable with a regular formula $M \Imp Z$
such that if $\VH$ is any Heyting algebra and $I$ is an interpretation in $\VH$
with $V_I(M) = \aT$, then $V_I(A) \le V_I(Z)$.
\end{Lemma}
\Proof
First assume $A$ is atomic.
If $A$ is $\lT$, let $Z$ be any variable and let $M \SDef Z$.
If $A$ is $\lF$, take $M$ and $Z$ to be distinct variables.
If $A$ is a variable, take $M$ to be some other variable and take $Z$
to be $A$.
In all three cases, $A$ and $M \Imp Z$ are either both provable
or both unprovable and hence they are equiprovable.

Now assume $A$ is not atomic.
By Lemma~\ref{lma:red}, we may assume $A$ is reduced.
If we choose some variable $Z$ that does not occur in $A$,
Then it is easy to see that $A \Eqp (A \Imp Z) \Imp Z$ (for the right-to-left
direction, use the substitution lemma to substitute $A$ for $Z$).
Our plan is to replace $K \SDef A\Imp Z$ by a basic context by ``unnesting''
all its non-atomic subformulas.
Assume $K$ contains $k$ non-atomic subformulas.
Starting with $K \equiv A_1 \equiv B_1 \circ_1 C_1$, enumerate
the $k$ non-atomic sub-formulas, $A_1 \equiv B_1 \circ_1 C_1, \ldots, A_k \equiv B_k \circ_k C_k$.
Choose fresh variables $P_i$, $i = 1, \ldots k$.
Define atomic formulas, $G_i$, $H_i$, for $i = 1, \ldots, k$ as follows:
$G_i$ is $B_i$ if $B_i$ is atomic and is $P_j$
if $B_i$ is the $j$-th non-atomic subformula;
$H_i$ is $C_i$ if $C_i$ is atomic and
is $P_j$ if $C_i$ is the $j$-th non-atomic subformula.
Now define formulas $L$ and $M$ as follows:
\begin{align*}
L &\SDef \bigwedge_{i=1}^{k} (P_i \Iff (G_i \circ_i H_i)) \\
M &\SDef P_1 \land L
\end{align*}

Recalling that $B \Iff C$ is just shorthand for $(B \Imp C) \And (C \Imp B)$,
and using the fact that $A$ and hence $K$ are reduced, we see that
$M$ is a basic context, so $M \Imp Z$ is regular.

We must show that $K \Imp Z \Eqp M \Imp Z$.
To see this, first assume $\IPLpr K \Imp Z$.
By induction on the size of the $A_i$, we have that
$\IPLpr L \Imp (P_i \Iff A_i)$, $i = 1, \ldots, k$.
Hence, as $\IPLpr M \Imp L$, $\IPLpr M \Imp (P_1 \Iff A_1)$,
i.e., $\IPLpr M \Imp (P_1 \Iff K)$.
As, clearly, $\IPLpr M \Imp P_1$, we have $\IPLpr M \Imp K$
and then, as $\IPLpr K \Imp Z$ by assumption, we have $\IPLpr M \Imp Z$.
Conversely, assume $\IPLpr M \Imp Z$. Using the substitution lemma,
we have also that $\IPLpr M[A_1/P_1, \ldots, A_k/P_k] \Imp Z$,
but $M[A_1/P_1, \ldots, A_k/P_k]$ is $K \land L'$
where $L' \equiv L[A_2/P_2, \ldots A_k/P_K]$ is a conjunction of
formulas of the form $A \Iff A$, hence
$\IPLpr M[A_1/P_1, \ldots, A_k/P_k] \Iff K$, and as
$\IPLpr M[A_1/P_1, \ldots, A_k/P_k] \Imp Z$ we have that $\IPLpr K \Imp Z$.

The claim about interpretations is clear for our choice of $M$ and $Z$
when $A$ is atomic.
In the case when $A$ is not atomic, construct $M$ and $Z$ as
described above and
assume $I$ is an interpretation such that $V_I(M) = \aT$.
Then for each $i = 1, \ldots, k$, we have
$V_I(P_i \Iff (G_i \circ_i H_i)) = \aT$, but this implies
that $V_I(P_i) = V_I(G_i \circ_i H_i)$ and hence,
(by induction on the size of the $A_i$) that $V_I(P_i) = V_I(A_i)$.
In particular, $V_I(P_1) = V_I(A_1)$ and since
we also have $V_I(P_1) = \aT$, we must have
$V_I(A_1) = \aT$. But by construction $A_1 \equiv A \Imp Z$,
so $V_I(A \Imp Z) = \aT$, which implies $V_I(A) \le V_I(Z)$.
\Done

We now state and prove three lemmas whose purpose will become clear
at their point of use in the proof of our main theorem, Theorem~\ref{thm:main}.

\begin{Lemma}\label{lma:simp-var}
If $B$ is a basic formula that is not of the form $P$ or $P \Imp Q \Or R$ and
$P$ occurs in $B$, then $\IPLpr P \And B \Iff P \And C$ where $C$ has fewer
connective occurrences than $B$ and is either a basic formula, an atom or a
basic context comprising a conjunction of two variables.
\end{Lemma}
\Proof Routine using the fact that $\IPL \Pr P
\And B \Iff P \And B[\lT/P]$ (which may be proved for arbitrary $B$ by
induction on the structure of $B$).
\Done

\begin{Lemma}\label{lma:simp-imp-imp}
If $\IPLpr K \And A \And (B \Imp C) \Imp B$, then
$
\IPL \Pr ((K \And ((A \Imp B) \Imp C)) \Imp D) \Iff (K \And C \Imp D).
$
\end{Lemma}
\Proof
$\Imp$: easy using $\IPLpr C \Imp ((A \Imp B) \Imp C)$. \\
$\Pmi$: the outline of a natural deduction proof is shown
in table~\ref{tab:simp-imp-imp}.
Here in step~\ref{labB}
we use $\IPLpr ((A \Imp B) \Imp C) \Imp (B \Imp C)$
to strengthen the antecedent of the implication.
\Done

\begin{table}
\begin{center}
\begin{align}
K \And A \And (B \Imp C) \Imp B &
		\mbox{\, [Given]} \label{Given}  \\
K \And C \Imp D &
		\mbox{\, [Assume]} \label{Assume} \\
K \And (B \Imp C) \Imp A \Imp B &
		\mbox{\, \eqref{Given}} \label{labA}\\
K \And ((A \Imp B) \Imp C) \Imp A \Imp B &
		\mbox{\, \eqref{labA}} \label{labB}\\
K \And ((A \Imp B) \Imp C) \Imp C &
		\mbox{\, \eqref{labB}} \label{labC}\\
K \And ((A \Imp B) \Imp C) \Imp D &
		\mbox{\, \eqref{labC} and \eqref{Assume}} \label{labD}\\
(K \And C \Imp D)  \Imp ((K \And ((A \Imp B) \Imp C)) \Imp D) &
		\mbox{\, \eqref{labD}, disch. \eqref{Assume}}
\end{align}
\caption{Outline natural deduction proof}
\label{tab:simp-imp-imp}
\end{center}

\end{table}

%\begin{figure}
%\begin{prooftree}
%   \[
%         \[
%            \[
%               \[
%                  \[
%                  \justifies
%                     K \And A \And (B \Imp C) \Imp B
%                  \using
%                     \mbox{Given}
%                  \]
%               \justifies
%                  K \And (B \Imp C) \Imp A \Imp B
%%               \using 
%%                  \mbox{\ldots Curry \ldots}
%               \]
%            \justifies
%               K \And ((A \Imp B) \Imp C) \Imp A \Imp B
%            \using
%                  \mbox{(*)}
%            \]
%         \justifies
%            K \And ((A \Imp B) \Imp C) \Imp C
%%         \using
%%            \mbox{\ldots $\Imp$-intro \ldots}
%         \]
%         [K \And C \Imp D]^{\alpha}
%      \justifies
%         K \And ((A \Imp B) \Imp C) \Imp D
%%      \using
%%         \mbox{\ldots $\Imp$-intro \ldots}
%      \]
%\justifies
%   (K \And C \Imp D)  \Imp ((K \And ((A \Imp B) \Imp C)) \Imp D)
%\using
%   \mbox{$\Imp$-intro}^{\alpha}
%\end{prooftree}
%\caption{Outline proof tree}
%\label{fig:simp-imp-imp}
%\end{figure}

\begin{Lemma}\label{lma:lift-counter-model}
Let $B$ be a basic formula that is not a variable and let
$I$ be an interpretation in a non-trivial Heyting
algebra $\VH$ such that $V_I(B) = \aT$.
Let $\alpha : H \to (H \Diff \{\aT\}) \cup \{\Star_H\}$ be as in the definition
of $\Gamma(\VH)$.
Define an interpretation $J$ in $\Gamma(\VH)$ by $J = \alpha \circ I$.

{\em(i)}
If $B$ does not have the form $(P \Imp E) \Imp R$ then $V_J(B) = \aT$.

{\em(ii)}
If $B$ has the form $(P \Imp E) \Imp R$, and if in addition $V_I(P) = V_I(E \Imp R) = \aT$ while $V_I(E) \neq \aT$,
then also $V_J(B) = \aT$.
\end{Lemma}
\Proof {\em(i)}:
This is easily checked for the case $P \Imp E$ and for the cases $P
\mathrel{\circ} Q \Imp R$ and $P \Imp Q \mathrel{\circ} R$ when ${\circ} \in
\{\And, \Or\}$.  In the remaining case $B \equiv P \Imp Q \Imp E$. As $B$ is
equivalent to $P \And Q \Imp E$, we have already covered the case when $E$ is a
variable, while if $E$ is $\lF$, $V_J(B) = \alpha(p) \Meet \alpha(q) \Res \aF$,
where $p = I(P)$ and $q = I(Q)$, but then, by inspection of the operation
tables, we have $\alpha(p) \Meet \alpha(q) = p \Meet q$ unless $p = q = \aT$,
but as $H$ is non-trivial and $V_I(B) = \aT$, the case $p = q = \aT$ cannot
arise.

{\em(ii)}:
we have $V_J(B) = (\alpha(p) \Res \alpha(e)) \Res \alpha(r)$, where
$p = V_I(P)$, $e = V_I(E)$ and $r = V_I(R)$. By assumption,
$p = \aT$ and $e \neq \aT$, so $\alpha(p) = \Star$ and $\alpha(e) = e$,
hence $\alpha(p) \Res \alpha(e) = \Star \Res e = e$,
so that $V_J(B) = e \Res \alpha(r)$ which is $e \Res \Star = \aT$,
if $r = \aT$, and is $e \Res r$ otherwise, in which case, as  we are given that $V_I(E \Imp R) = \aT$, we have $e \Res r = V_I(E \Imp R) = \aT$.
\Done

To state our main theorem, we define an interpretation $I$ to be a {\em strong
refutation} of a formula of the form $K \Imp C$, if $V_I(K) = \aT$
while $V_I(C) \neq \aT$.

\begin{Theorem}\label{thm:main}
Let $A \equiv K \Imp F$ be a regular formula (so that $F$ is either a variable or $\lF$), let
$K \equiv B_1 \And \ldots \And B_k$ display $K$ as a disjunction of basic
formulas and let $d = d(A)$ be the number of $B_i$ of the form
$(P \Imp E) \Imp R$.
Either $\IPL \Pr A$ or $A$ has a strong refutation in $\Jas_d$.
\end{Theorem}
\Proof
The proof is by induction on the sum $s(A)
= c(A) + d(A) + v(A)$, where $c(A)$ is the number of connective occurrences in
$K$, $d(A)$ is as in the statement of the theorem and $v(A)$ is the number of
conjuncts of $K$ comprising a single variable.

Case {\Lab{i}}: $v(A) = d(A) =  0$: in this case, the interpretation
in $\Jas_0 = \BB$ that maps every variable to $\aF$ is easily seen
to be a strong refutation of $A$ (which is therefore unprovable,
by the soundness of $\IPL$).

Case {\Lab{ii}}: $v(A) > 0$: in this
case at least one $B_i$ is a variable.
If all the $B_i$ are variables and if $B_i \not\equiv F$ for any $i$, then $A$ has strong refutation
such that $I(B_i) = \aT$, $i = 1, \ldots, k$ and $V_I(F) = \aF$.
Otherwise, rearranging the $B_i$ if necessary, we may assume that
$K \equiv P \And L$ where $P$ is a variable and
$L \equiv B_2 \And \ldots \And B_k$.
If $P \equiv F$, we are done: $F \And L \Imp F$ is provable.
If $P \not\equiv F$ and $P$ does not occur in $L$, then it is easy
to see that $A \Eqp A'$ where $A' \SDef L \Imp F$.
As $s(A') < s(A)$, by induction, if $\IPL \not\Pr L \Imp F$,
we can find a strong refutation $I$ of $L \Imp F$, but then, because $P$ does
not occur in $L \Imp F$, by adjusting $I$ if necessary to map $P$ to $\aT$ we obtain a strong refutation of $A$.
If $P$ occurs in $L$, let us rearrange the $B_i$ again so that
$K \equiv P \And B \And M$ where $M \equiv  B_3, \ldots, B_k$
and $P$ occurs in $B$.
If $B$ does not have the form $P \Imp Q \Or R$, then,
by Lemma~\ref{lma:simp-var}, we may replace $P \And B$ by an equivalent
formula $P \And C$ where $C$ is either a basic formula, an atom or a basic
context comprising a conjunction of two variables
and contains fewer connectives then $B$. If $C$ is $\lF$, $A$ is provable and we are done.
Otherwise, we may replace $A$ by
the equivalent regular formula $A' \SDef P \And C \And M \Imp F$ (or $P \And M \Imp F$, if $C$ is $\lT$)
and we are done by induction, since $s(A') < s(A)$.
If $B$ has the form $P \Imp Q \Or R$, then
$\IPL \Pr P \And B \And M \Iff K' \Or K''$
where $K' \SDef P \And Q \And M$ and $K'' \SDef P \And R \And M$,
and hence $\IPL \Pr A \Iff A' \And A''$
where $A' \SDef K' \Imp F$ and $A'' \SDef K'' \Imp F$.
If $A$ is not provable, then one of $A'$ and $A''$ is not provable,
in which case, as $s(A') < s(A)$ and $s(A'') < s(A)$, by induction
we have a strong refutation in $\Jas_d$ of either $A'$ or $A''$ and this
will also strongly refute $A$.

Case {\Lab{iii}}: $v(A) = 0$ and $d = d(A) > 0$:
Let $X = \{j_1, \ldots, j_d\}$ be the set of $i$ such that $B_i$ has the form
$(P \Imp E) \Imp R$.  For each $i \in X$, let $K_i \SDef B_1 \And \ldots \And
B_{i-1} \And B_{i+1} \And \ldots \And B_k$ and let $P_i$, $E_i$ and $R_i$ be
such that $B_i \equiv (P_i \Imp E_i) \Imp R_i$.
We now have two subcases depending on the provability of the formulas
$C_i \SDef K_i \And P_i \And (E_i \Imp R_i) \Imp E_i$:

Subcase {\Lab{iii}}{\em(a)}:
for some $i \in X$, $\IPLpr C_i$:
By Lemma~\ref{lma:simp-imp-imp},
$A$, which is equivalent to $K_i \And ((P_i \Imp E_i) \Imp R_i) \Imp F$,
is equivalent to $A' \SDef K_i \And R_i \Imp F$.  As $s(A') < s(A)$,
we are done by induction.

Subcase {\Lab{iii}}{\em(b)}:
for every $i \in X$, $\IPL \not\Pr C_i$:
By induction, as $s(C_i) < s(A)$ and $d(C_i) = d - 1$,
for each $i \in X$ there is an interpretation
$I_i$ in $\Jas_{d-1}$ that strongly refutes $C_i$, i.e., $K_i \And P_i \And (E_i
\Imp R_i) \Imp E_i$.  Now define an interpretation $I$ in $\Jas_{d-1}^d$, by
$I(U) = (I_{j_1}(U), \ldots, I_{j_d}(U))$.  Then $V_I(B_i) = \aT$ for
$i = 1, \ldots, k$ (because,
for $i \in X$, $V_{I_i}(P_i) = V_{I_i}(E_i \Imp R_i) = \aT$
and $B_i \equiv (P_i \Imp E_i) \Imp R_i$).
But then applying Lemma~\ref{lma:lift-counter-model} to $I$ gives
us an intepretation $J$ in $\Jas_d = \Gamma(\Jas_{d-1}^d)$ that strongly
refutes $A$.
\Done

\begin{Corollary}\label{cor:main-iff}
Let $A \equiv K \Imp F$ be a regular formula and let $d$ be the number of conjuncts
of $K$ of the form $(P \Imp E) \Imp R$. Then $\IPL \Pr A$ iff $\Jas_d \models A$.
\end{Corollary}
\Proof Immediate from the theorem given the soundness of $\IPL$ for the Heyting algebra semantics.
\Done

\begin{Corollary}\label{cor:complete}
$\IPL$ is complete for the Heyting algebra semantics.
\end{Corollary}
\Proof Assume $\models A$. We have to show that $\IPLpr A$.
Consider the regular formula $A' \equiv M \Imp Z$ such that $A \Eqp A'$
given by Lemma~\ref{lma:reg}.
If $\IPL \not\Pr A$, then $\IPL \not\Pr A'$, whence by the theorem,
$A'$ has a strong refutation in $\Jas_k$ for some $k$, i.e., an interpretation
$I$ in $\Jas_k$ such that $V_I(M) = \aT$, but $V_I(Z) < \aT$.
But then Lemma~\ref{lma:reg} gives us that $V_I(A) \le V_I(Z) < \aT$,
so $I \not\models A$ contradicting our assumption that $\models A$.
\Done

\begin{Corollary}\label{cor:FMP}
$\IPL$ has the finite model property.
\end{Corollary}
\Proof From the theorem and soundness we know that a refutable
regular formula has a refutation in a finite model. Argue as in the proof
of Corollary~\ref{cor:complete} to reduce the general case to the
case of regular formulas.
\Done

If $\VH_0, \VH_1, \ldots$ is a sequence of Heyting algebras, let
us define $\bigodot_k \VH_k$ to be the subalgebra of $\prod_k \VH_k$
comprising sequences $(p_0, p_1, \ldots)$ such that for all sufficiently
large $k$, the $p_k$ are either all $\aF$ or all $\aT$.
Our final corollary shows that there is countably infinite Heyting algebra
$\Jas$, such that for any formula $\phi$, $\Jas \models \phi$ iff
$\IPL \vdash \phi$.

\begin{Corollary}\label{cor:J-complete}
For any formula $A$, $\IPLpr A$ iff $\bigodot_k \Jas_k \models A$.
\end{Corollary}
\Proof
The left-to-right direction is just the soundness of $\IPL$ for Heyting algebras.
For the right-to-left direction argue as in the proof of Corollary~\ref{cor:complete} and note that a refutation in $\Jas_d$ gives a refutation in the
subalgebra of $\bigodot_k \Jas_k$ comprising the sequences $(p_0, p_1, \ldots)$ such that
$p_i$ is constant for $i > d$.
\Done

The statement of Theorem~\ref{thm:main}
leads to a decision procedure for $\IPL$ that involves 
a search through all interpretations of a formula in one of the $\Jas_d$
for a certain $d$. As Rose \cite{Rose53} observes, the size of the $\Jas_k$
grows very rapidly with $k$, so this decision procedure is impractical.
However, the proof of the theorem leads to a much better algorithm:
given any formula $A$, we first apply the algorithm of Lemma~\ref{lma:reg} if necessary
to convert $A$ into an equiprovable regular formula and then follow the case analysis
of the proof of the theorem: if we are in Case~{\Lab{i}}, $A$ is unprovable and we are done;
if we are in Case~{\Lab{ii}}, the proof shows us how to produce one or two simpler
formulas whose conjunction is equivalent to $A$ and we may proceed recursively
to decide these formulas; if we are in Case~{\Lab{iii}}, we can derive the formulas
$C_i$ described in the proof and decide them recursively; if any $C_i$ is provable,
we are in Subcase~{\Lab{iii}}{\em(a)} and we may replace $A$ by an equivalent
and simpler formula that we can decide recursively; if no $C_i$ is provable,
we are in Subcase~{\Lab{iii}}{\em(b)} and $A$ is unprovable.
If $A$ is unprovable, then the proof of the theorem yields
an explicit refutation in one of the $\Jas_k$.
In the appendix, we show some example calculations using this decision procedure.
We make no claim that the decision procedure is practical on large examples:
its time complexity involves a factor $d!$, where $d$ is bounded below by the
number of implications in the input formula.

\Jaskowski's construction was used by Tarski to show the completeness of
intuitionistic propositional logic for its topological interpretation
\cite{Tarski38}. One imagines that the details of the proof that {\Jaskowski} sketched in \cite{Jaskowski36} were well known to Polish
logicians in the 1930s, but sadly the details have been lost: by the 1950s,
Kleene's student Gene F. Rose had to reinvent a proof.
The proof of Theorem~\ref{thm:main} given here and, in particular, its use of
Lemma~\ref{lma:simp-imp-imp} is largely due to Rose \cite{Rose52,Rose53}.
Rose's analogue of our notion of basic formula admits only 6 forms:
$P$, $\Not P$, $P \Imp Q$, $P \Imp Q \Or R$, $P \And Q \Imp R$ and
$(P \Imp Q) \Imp R)$. To prove his analogue of our Lemma~\ref{lma:reg} involves
a lengthy case analysis, whereas our more liberal notion of basic formula
admits the simpler and more intuitive proof given here. As far as I know,
the observations that Theorem~\ref{thm:main} leads to an alternative proof
of the completeness of $\IPL$ and that its proof leads to a syntax-driven
decision procedure for $\IPL$ are new.

\bibliographystyle{asl}
\bibliography{bookspapers}

%\newpage
\appendix
\section*{Appendix: examples of the decision procedure}
Throughout the examples ``Case'' and ``Subcase'' refer to the proof of
Theorem~\ref{thm:main}.
We use the following tabular format for the regular formulas $B_1 \And \ldots \And B_k \Imp F$
that occur as the goals we are trying to decide:
$$
\Reg{B_1, \ldots, B_k}{F}
$$

\subsection*{Example 1: $A \SDef (P \Or Q) \And \Not Q \Imp P$}
Noting that $A$ already has the form $B \Imp Q$, we can skip the first step in the algorithm of Lemma~\ref{lma:reg}
and simply ``unnest'' $B$. Listing the subformulas of $(P \Or_2 Q) \And_1 \Not_3 Q$ as shown by the subscripts,
our initial goal is:
$$
\Reg{P_1, P_1 \Iff P_2 \And P_3, P_2 \Iff P \Or Q, P_3 \Iff \Not Q}{P}
$$
We are in Case~\Lab{ii} and we replace the occurrence of $P_1$ in $P_1 \Iff
P_2 \And P_3$ by $\lT$ and simplify giving;
$$
\Reg{P_1, P_2, P_3, P_2 \Iff P \Or Q, P_3 \Iff \Not Q}{P}
$$
We are again in Case~\Lab{ii}, but now $P_2$ appears in a subformula of the form $P_2 \Imp P \Or Q$ and replacing
$P_2$ by $\lT$ in that formula gives us two subgoals:
$$
\Reg{P_1, P_2, P_3, P, P_3 \Iff \Not Q}{P} \quad
\Reg{P_1, P_2, P_3, Q, P_3 \Iff \Not Q}{P}
$$
Both subgoals are in Case~\Lab{ii}. In the first, the succedent of the goal appears in the antecedent
while in the second, replacing first $P_3$ and then $Q$ by $\lT$ in $P_3 \Iff \Not Q$ and simplifying gives the antecedent $\lF$.
So both subgoals and hence also our original formula are provable.

\section*{Example 2: Peirce's law: $A \SDef ((P \Imp Q) \Imp P) \Imp P$}
$A$ is already regular, so we take it as our initial goal:
$$
\Reg{(P \Imp Q) \Imp P}{P}
$$
We are in Case~\Lab{iii} and our next step is to decide the goal:
$$
\Reg{P, Q \Imp P}{Q}
$$
This is in Case~\Lab{ii} and replacing $P$ by $\lT$ in $Q \Imp P$ and simplifying leads to
$$
\Reg{P}{Q}
$$
This is again in Case~\Lab{ii} and is refuted by the interpretation $\{P \mapsto \aT, Q \mapsto \aF\}$.
Following Lemma~\ref{lma:lift-counter-model}, this lifts to the refutation
$\{P \mapsto *, Q \mapsto \aF\}$ of Peirce's law in $\Jas_1 = \BB \cup \{*\}$.

\section*{Example 3: prelinearity: $A \SDef (P \Imp Q) \Or (Q \Imp P)$}
Following the first part of Lemma~\ref{lma:reg}, we replace $A$
by the equiprovable formula $(A \Imp Z) \Imp Z$ and list its
subformulas as indicated by the subscripts in
$((P \Imp_3 Q) \Or_2 (Q \Imp_4 P) \Imp_1 Z) \Imp Z$. This
gives us the following initial goal:
$$
\Reg{P_1, P_1 \Iff P_2 \Imp Z, P_2 \Iff P_3 \Or P_4,
   P_3 \Iff P \Imp Q, P_4 \Iff (Q \Imp P)}{Z}
$$
This is in Case~\Lab{ii} and replacing $P_1$ by $\lT$ in $P_1 \Iff P_2 \Imp Z$
and simplifying we get:
$$
\Reg{P_1, P_2 \Imp Z, P_2 \Iff P_3 \Or P_4,
   P_3 \Iff P \Imp Q, P_4 \Iff (Q \Imp P)}{Z}
$$
This is now in Case~\Lab{iii} with $d = 2$. This leads to two subgoals:
$$
\begin{array}{rl}
C_1{:} & \Reg{%
   \begin{array}{@{}c@{}}
      P_1, P_2 \Imp Z, P_2 \Iff P_3 \Or P_4, P_3 \Imp P \Imp Q,\\
      P_4 \Iff (Q \Imp P), P, Q \Imp P_3
   \end{array}}{Q} \\\ \\
C_2{:} & \Reg{%
   \begin{array}{@{}c@{}}
      P_1, P_2 \Imp Z, P_2 \Iff P_3 \Or P_4, P_3 \Iff P \Imp Q,\\
      P_4 \Imp (Q \Imp P), Q, P \Imp P_4
   \end{array}}{P} \\\ \\
\end{array}
$$
Either continuing to follow Theorem~\ref{thm:main} or by inspection, we find
the following strong refutations of these subgoals in $\BB$.
$$
\begin{array}{rl}
C_1{:} & (\{P, P_1, P_2, P_4, Z\} \times \{\aT\}) \cup (\{Q, P_3\} \times \{\aF\}) \\
C_2{:} & (\{Q, P_1, P_2, P_3, Z\} \times \{\aT\}) \cup (\{P, P_4\} \times \{\aF\})
\end{array}
$$
Combining these we should obtain a refutation $I = \{P \mapsto (\aT, \aF), Q \mapsto
(\aF, \aT)\}$ of $A$ in $\Gamma(\BB^2) \subseteq
\Jas_2$. And, indeed, in $\Gamma(\BB^2)$ we have:
\begin{align*}
((\aT, \aF) \Res (\aF, \aT)) \Join ((\aF, \aT) \Res (\aT, \aF))
   & = (\aF, \aT) \Join (\aT, \aF) \\
   & = \alpha((\aF, \aT) \Join_{\BB^2} (\aT, \aF)) \\
   & = \alpha((\aT, \aT)) = * \neq \aT.
\end{align*}

\end{document}